\documentclass[12pt]{article}
\usepackage{amsfonts}
\usepackage{mathrsfs}
\usepackage{amsmath}
\usepackage{latexsym}
\usepackage{amssymb}
\usepackage{amsthm}
\usepackage{amscd}

\title{Classification of  Quiver Hopf Algebras and Pointed Hopf Algebras of  Type One}
\author{\small
Shouchuan Zhang $^{a, c}$,  Hui-Xiang Chen $^b$,  Yao-Zhong Zhang $^c$,   \ \  \ \\
\small $a$. Department  of Mathematics,
Hunan University\\  \small   Changsha  410082, \ P.R. China \\
\small $b$. Department of Mathematics, Yangzhou University \\
 \small   Yangzhou 225002, P.R. China\\
\small $c$. School of Mathematics and Physics, The University of Queensland\\
\small Brisbane 4072, Australia.
}
\date{}

\begin{document}
\newtheorem{Proposition}{Proposition}[section]
\newtheorem{Theorem}{Theorem}
\newtheorem{Definition}[Proposition]{Definition}
\newtheorem{Corollary}[Proposition]{Corollary}
\newtheorem{Lemma}[Proposition]{Lemma}
\newtheorem{Example}[Proposition]{Example}

\maketitle 

\numberwithin{equation}{section}

\date{}

\begin{abstract} Quiver Hopf algebras are classified by means
of ramification systems with irreducible representations. This leads
to the classification of  Nichols algebras over group algebras and
pointed Hopf algebras of  type one.

\vskip0.5cm
{\small 

}
\end {abstract}

\section{\bf Introduction}

Hopf algebras have important applications in mathematics and
mathematical physics. Indeed, quasi-triangular Hopf algebras give rise to braided  tensor categories 
and Yang-Baxter equations of important applications in integrable systems and
statistical mechanics (see. e.g. \cite{Majid90} and references therein).
Semi-simple Hopf algebras and non-semisimple Hopf
algebras  are related to  conformal field theories (see e.g. \cite{Ga}).

Classification of Hopf algebras is a main objective in the research
of Hopf algebras. So far, many important results have been obtained
in the classification of finite dimensional pointed Hopf algebras
(see e.g. \cite { AS98b, AS02, AS05, He04, AZ07}). Classification of PM quiver Hopf algebras
was completed by means of ramification system with characters 
in \cite {ZZC04}. Classification of ramification systems with characters
over the symmetric group $\mathbb S_n$ with $n\not = 6$ was obtained in
\cite {ZWW08}. Irreducible Hopf bimodules over a finite group were
described in \cite {DPR}. They correspond to pairs $(\mathcal O_s,
\rho )$,  where $\mathcal O_s$ is a conjugacy class containing $s$
in G and $\rho$ is an irreducible representation of the centralizer $G^s$  in G. 

In this paper,   quiver Hopf algebras, Nichols algebras over group
algebras and pointed Hopf algebras of  type one are classified by
means of ramification system with irreducible representations ({\rm
RSR} in short). As examples we classify RSRs over
symmetric group $\mathbb S_n$ with $n \not=6.$

Quivers \cite {CR02, CR97, ZZ07, OZ04} and tensor algebras of Hopf bimodules \cite{Ni78, Wo89} 
have widely been applied in representation theory,  Hopf algebras and quantum groups.
We use quivers to describe Yetter-Drinfeld $kG$-modules, $kG$-Hopf bimodules, Nichols algebras in braided
tensor category $^{kG}_{kG} {\mathcal YD}$, pointed  Hopf algebras of type one and   quiver Hopf algebras.

\section*{\bf Preliminaries}\label{s0}

Throughout this paper we assume that $G$ is a finite group and $k$ is a field. 

 Let $\widehat{{ G}}$ denote the set
of all isomorphism classes of irreducible representations of group
$G$ and  $Z_{s}$  the centralizer of $s$ in $G$. For $h \in G$ and
an isomorphism $\phi $ from $G$ to $G'$, define a map $\phi _h$ from
$G$ to $G'$ by sending $x$ to $\phi (h^{-1}xh)$ for any $x\in G$. ${\rm deg} (\rho)$
denotes the dimension of a representation space of representation  $\rho$.

Let ${\mathbb N}$ and ${\mathbb Z}$ denote the sets of  all non-negative integers and
all integers, respectively. For a set $X$,
we denote by $|X|$ the  number of elements in $X$.  If $X = \oplus _{i\in I} X_{(i)}$ as vector
spaces, then we denote by $\iota _i$ the natural injection from
$X_{(i)}$ to $X$ and by $\pi _i$ the corresponding projection from $X$
to $X_{(i)}$. We will use $\mu$ to denote the multiplication  of an
algebra and use $\Delta$ to denote the comultiplication of a
coalgebra. For a (left or right) module and a (left or right)
comodule, denote by $\alpha ^-$, $\alpha ^+$, $\delta ^-$ and
$\delta ^+$ the left module, right module, left comodule and right
comodule structure maps, respectively. Sweedler's sigma
notations for coalgebras and comodules are $\Delta (x) = \sum
x_{(1)}\otimes x_{(2)}$, $\delta ^- (x)= \sum x_{(-1)} \otimes
x_{(0)}$, $\delta ^+ (x)= \sum x_{(0)} \otimes x_{(1)}$.

A quiver $Q=(Q_0,Q_1,s,t)$ is an oriented graph, where  $Q_0$ and
$Q_1$ are the sets of vertices and arrows, respectively; $s$ and $t$
are two maps from  $Q_1$ to $Q_0$. For any arrow $a \in Q_1$, $s(a)$
and $t(a)$ are called its start vertex and end vertex, respectively,
and $a$ is called an arrow from $s(a)$ to $t(a)$. For any $n\geq 0$,
an $n$-path or a path of length $n$ in the quiver $Q$ is an ordered
sequence of arrows $p=a_na_{n-1}\cdots a_1$ with $t(a_i)=s(a_{i+1})$
for all $1\leq i\leq n-1$. Note that a 0-path is exactly a vertex
and a 1-path is exactly an arrow. In this case, we define
$s(p)=s(a_1)$, the start vertex of $p$, and $t(p)=t(a_n)$, the end
vertex of $p$. For a 0-path $x$, we have $s(x)=t(x)=x$. Let $Q_n$ be
the set of $n$-paths. Let $^yQ_n^x$ denote the set of all $n$-paths
from $x$ to $y$, $x, y\in Q_0$. That is, $^yQ_n^x=\{p\in Q_n\mid
s(p)=x, t(p)=y\}$.

A quiver $Q$ is {\it finite} if $Q_0$ and $Q_1$ are finite sets. A
quiver $Q$ is {\it locally finite} if $^yQ_1^x$ is a finite set for
any $x, y\in Q_0$.

Let $G$ be a group. Let ${\mathcal K}(G)$ denote the set of
conjugacy classes in $G$. A formal sum $r=\sum_{C\in {\mathcal
K}(G)}r_CC$  of conjugacy classes of $G$  with cardinal number
coefficients is called a {\it ramification} (or {\it ramification
data} ) of $G$, i.e.  for any $C\in{\mathcal K}(G)$, \  $r_C$ is a
cardinal number. In particular, a formal sum $r=\sum_{C\in {\mathcal
K}(G)}r_CC$  of conjugacy classes of $G$ with non-negative integer
coefficients is a ramification of $G$.

 For any ramification $r$ and a $C \in {\mathcal K}(G)$, without loss of generality,
we can choose a set $I_C(r)$ such that its cardinal number is $r_C$.
 Let ${\mathcal K}_r(G):=\{C\in{\mathcal
K}(G)\mid r_C\not=0\}=\{C\in{\mathcal K}(G)\mid
I_C(r)\not=\emptyset\}$.  If there exists a ramification $r$ of $G$
such that the cardinal number of $^yQ_1^x$ is equal to $r_C$ for any
$x, y\in G$ with $x^{-1}y \in C\in {\mathcal K}(G)$, then $Q$ is
called a {\it Hopf quiver with respect to the ramification data
$r$}. In this case, there is a bijection from $I_C(r)$ to $^yQ_1^x$,
and hence we write  ${\ }^yQ_1^x=\{a_{y,x}^{(i)}\mid i\in I_C(r)\}$
for any $x, y\in G$ with $x^{-1}y \in C\in {\mathcal K}(G)$. Denote
by $ (Q, G, r)$ the Hopf quiver of $G$ with respect to $r$.

If  $\phi: A\rightarrow A'$ is  an algebra homomorphism and  $(M,
\alpha ^-)$ is a left $A'$-module, then $M$ becomes a left
$A$-module with the $A$-action given by $a \cdot x =\phi (a) \cdot x
$ for any $a\in A$, $x\in M$, called a pullback $A$-module through
$\phi$, written as  $_{\phi}M$.  Dually, if  $\phi: C\rightarrow C'$
is a coalgebra homomorphism and  $(M, \delta ^- )$ is  a left
$C$-comodule, then $M$ is a left $C'$-comodule with the
$C'$-comodule structure given by $ {\delta'}^-:=(\phi\otimes{\rm
id})\delta^-$, called  a push-out $C'$-comodule through $\phi$,
written as  $^{\phi}M$.

If $B$ is a Hopf algebra and $M$ is a $B$-Hopf bimodule, then  we
say that $(B, M)$ is a  Hopf bimodule. For any two Hopf bimodules
$(B,M)$ and $(B', M')$, if $\phi$ is a Hopf algebra homomorphism
from $B$ to $B'$  and $\psi$ is simultaneously a $B$-bimodule homomorphism from
$M$ to $_\phi M'{}_\phi $ and a $B'$-bicomodule homomorphism from
$^\phi M ^\phi $ to $M'$, then $(\phi, \psi)$ is called a pull-push
Hopf bimodule homomorphism. If $\psi$ is a bijection, then we say that $(\phi, \psi)$  is a
a pull-push Hopf bimodule isomorphism, written as $(B, M) \cong (B', M')$ as  pull-push
Hopf bimodules. In particular, if $B=B'$ we also write $ M \cong M'$ as  pull-push
$B$-Hopf bimodules, in short.  Similarly, we say that $(B, M)$ and $(B,
X)$ are a Yetter-Drinfeld module and a Yetter-Drinfeld Hopf algebra, respectively, if $M$ is a
Yetter-Drinfeld $B$-module and $X$ is a braided Hopf algebra in
Yetter-Drinfeld category $^B_B {\mathcal YD}$.  For any two
Yetter-Drinfeld modules $(B,M)$ and $(B', M')$, if $\phi$ is a Hopf algebra
homomorphism from $B$ to $B'$,  and $\psi$ is simultaneously a left
$B$-module homomorphism from  $M$ to $_\phi M' $ and a left
$B'$-comodule homomorphism from $^\phi M  $ to $M'$, then $(\phi,
\psi)$ is called a pull-push Yetter-Drinfeld module homomorphism.  For any
two Yetter-Drinfeld Hopf algebras $(B,X)$ and $(B', X')$, if $\phi$ is a
Hopf algebra homomorphism from $B$ to $B'$,   $\psi$ is
simultaneously a left $B$-module homomorphism from  $X$ to $_\phi X'
$ and a left $B'$-comodule homomorphism from $^\phi X  $ to $X'$,
meantime, $\psi$ also is algebra and coalgebra homomorphism from $X$
to $X'$, then $(\phi, \psi)$ is called a pull-push Yetter-Drinfeld Hopf
algebra homomorphism (see the remark after Theorem 4 in \cite{ZZC04}).

Let $A$ be an algebra and $M$ be an $A$-bimodule. Then the tensor
algebra $T_A(M)$ of $M$ over $A$ is a graded algebra with
$T_A(M)_{(0)}=A$, $T_A(M)_{(1)}=M$ and $T_A(M)_{(n)}=\otimes^n_AM$ for $n>1$.
That is, $T_A(M)=A\oplus(\bigoplus_{n>0}\otimes^n_AM)$ (see
\cite{Ni78}). Let $D$ be another algebra. If $h$ is an algebra map
from $A$ to $D$ and $f$ is an $A$-bimodule map from $M$ to $D=\
_hD_h$, then by the universal property of $T_A(M)$ (see Proposition 1.4.1 in \cite
{Ni78}) there is a unique algebra map
$T_A(h,f): T_A(M)\rightarrow D$ such that $T_A(h,f)\iota_0=h$ and
$T_A(h,f)\iota_1=f$. One can easily see
that $T_A ( h, f ) = h + \sum _{n>0} \mu ^{n-1}T_n (f )$. For the details, the reader is
directed to Section 1.4 of \cite{Ni78} or \cite {ZZC04}.
Dually, let $C$ be a coalgebra and let $M$ be a $C$-bicomodule.
Then the cotensor coalgebra $T_C^c(M)$ of $M$ over $C$ is a graded
coalgebra with $T_C^c(M)_{(0)}=C$, $T_C^c(M)_{(1)}=M$ and
$T_C^c(M)_{(n)}=\Box^n_CM$ for $n>1$. That is,
$T_C^c(M)=C\oplus(\bigoplus_{n>0}\Box^n_CM)$ (see \cite{Ni78} or \cite {ZZC04} ). If $B$ is a Hopf algebra
and $M$ is a $B$-Hopf bimodule, then both $T_B(M)$ and $T_B^c (M)$ are graded Hopf algebras. Furthermore,
the subalgebra generated by $H$ and $M$ in $T_B^c(M)$, written as $H[M]$,   is a Hopf subalgebra of
$T_B^c(M)$ and $H[M] $ is called a Hopf algebra of type one.

\section {\bf Classification of Quiver Hopf Algebras }\label{s1}

\begin{Definition}\label {1.1}
$(G, r, \overrightarrow \rho, u)$ is called a ramification system
with irreducible representations  (or {\rm RSR } in short), if $r$
is a ramification of $G$; $u$ is a map from ${\mathcal K}(G)$ to $G$
with $u(C)\in C$ for any $C\in {\mathcal K}(G)$;  $I_{C} (r, u)$ and
$J_C(i)$ are  sets   with $ \mid\!J_{C} (i )\!\mid $  = $ {\rm deg}
(\rho _C^{(i)})$ and $I_C(r) = \{(i, j) \mid i \in I_C(r, u), j\in
J_C(i)\}$  for any $C\in {\mathcal K}_r(G)$, $i\in I_C(r, u)$;
$\overrightarrow \rho=\{\rho_C^{(i)} \}_ { i\in I_C(r,u ),
C\in{\mathcal K}_r(G)} \ \in \prod _ { C\in{\mathcal K}_r(G)}
(\widehat{    {Z_{u(C)}}}) ^{\mid I_{C} (r, u) \mid }$ with
$\rho_C^{(i)} \in \widehat{ {Z_{u(C)}}} $ for any
 $ i \in I_C(r, u), C\in {\mathcal K}_r(G)$.
${\rm RSR} (G, r, \overrightarrow \rho, u)$ and ${\rm RSR} (G', r',
\overrightarrow {\rho '}, u')$ are said to be {\it isomorphic} if
the following conditions are satisfied:

$\bullet$ There exists a group isomorphism $\phi: G\rightarrow G'$.

$\bullet$ For any $C\in{\mathcal K}(G)$, there exists an element
$h_C\in G$ such that $\phi(h_C^{-1}u(C)h_C)=u'(\phi(C))$.

$\bullet$ For any $C\in{\mathcal K}_r(G)$, there exists a bijective
map $\phi_C : I_C(r, u) \rightarrow I_{\phi (C)}(r', u')$ such that
$\rho_C^{(i)} \cong\rho'{}_{\phi(C)}^{(\phi_C(i))} \phi _{h_c}$ as representations of
 $kZ_{u(C)}$ for all $i\in I_C(r, u)$, where $\phi
_{h_C} (h) = \phi (h^{-1}_C h h_C)$ for any $h \in G$.
\end {Definition}

\noindent {\bf Remark.} Assume that $G=G'$, $r=r'$,  $u(C) = u'(C)$ and $I_C(r, u)=
I_C(r', u')$ for any
$C \in {\mathcal K}_r(G)$. If there is a permutation $\phi_C$ on $I_C(r,
u)$ for any $C\in{\mathcal K}_r(G)$ such that
$\rho'{_C^{(\phi_C(i))}}\cong\rho_C^{(i)}$ for all $i\in I_C(r, u)$, then obviously
 ${\rm RSR} (G, r, \overrightarrow \rho, u)\cong {\rm RSR} (G, r, \overrightarrow{\rho '}, u)$.

\begin {Example}\label{1.4} Assume that  $k$ is a complex field  and $G = \mathbb
{S}_3$, then there are 3 elements in ${\mathcal K}(G)$, which are
$\{(1)\}$, $\{(12),(13),(23)\}$, $\{(123),(132)\}$,  and there are 3
non-isomorphic irreducible representations, namely two   1
dimensional irreducible representations,   $\epsilon$ and  $ \rm
{sgn} $, and one 2 dimensional irreducible representation $\rho$.
Obviously $Z_{u(\{1\})}=G$.
The set $\{{\rm RSR} (G, r, \overrightarrow \rho, u) $ $ \mid \overrightarrow
\rho = \rho, (sgn, sgn), (\epsilon, sgn ),  (\epsilon, \epsilon )
\}$ gives all the representatives of isomorphic classes of  ${\rm RSR}$s
with $r = r_Cr$ and $C= \{(1)\}$. Furthermore, when $\overrightarrow
\rho = \rho$, we can set $I_C(r, u) = \{ 1 \}$ and $J_C(1) = \{1
\}$. In this case $\rho _C^{(1)} =\rho$. When $\overrightarrow \rho
= (\epsilon, sgn )  $, we can set $I_C(r, u) = \{ 1, 2 \}$ and
$J_C(1) = J_C(2)= \{1 \}$. In this case $\rho _C^{(1)} =\epsilon$,
$\rho _C^{(2)} =sgn$.

\end {Example}

Let \begin {eqnarray} \label {e0.1}
G &=&\bigcup_{\theta\in\Theta_C}Z_{u(C)}g_{\theta},
\end {eqnarray}
where $\Theta_C$ is an index set, be  a coset decomposition of $Z_{u(C)}$ in $G$.
It is easy to check that
$|\Theta_C|=|C|$. We always assume that the representative element
of the coset $Z_{u(C)}$ is the identity $1$ of $G$.
 For any $h\in G$ and
$\theta\in\Theta_C$, there exist unique $h'\in Z_{u(C)}$ and
$\theta'\in\Theta_C$ such that $g_{\theta}h = h'g_{\theta'}$. Let
$\zeta_{\theta}(h)=h'$.  Then we have
\begin {eqnarray} \label {e0.3} g_{\theta}h&=&\zeta_{\theta}(h)g_{\theta'}.
\end {eqnarray}

 Let $Q=(G, Q_1, s, t)$ be a quiver of a group $G$. Then $kQ_1$
becomes  a $kG$-bicomodule  under  the natural comodule structures:
\begin{eqnarray}\label{arcom}
\delta^-(a)=t(a)\otimes a,\ \ \delta^+(a)=a\otimes s(a),\ \ a\in
Q_1,
\end{eqnarray} called an {\it arrow comodule}, written as $kQ_1 ^c$. In this case,
the path coalgebra $kQ^c$ is exactly isomorphic to the cotensor
coalgebra $T^c_{kG}(kQ_1^c)$  over $kG$ in a natural way (see
\cite{CR97} and \cite{CR02}).  We will set $kQ^c=T^c_{kG}(kQ_1^c)$
in the following. Moreover, when $G$ is finite,  $kQ_1$ becomes a
$(kG)^*$-bimodule with the module structures defined by
\begin{eqnarray}\label{armod}
\mbox{\hspace{1cm}}p\cdot a:=\langle p, t(a)\rangle a,\ \ a\cdot
p:=\langle p, s(a)\rangle a,\ \ p\in(kG)^*, a\in Q_1,
\end {eqnarray}
written as $kQ_1^a$, called an {\it arrow module}. Therefore, we
have a tensor algebra $T_{(kG)^*}(kQ_1^a)$. Note that the tensor
algebra $T_{(kG)^*}(kQ_1^a)$ of $kQ_1^a$ over $(kG)^*$ is exactly
isomorphic to the path algebra $kQ^a$. We will set
$kQ^a=T_{(kG)^*}(kQ_1^a)$ in the following.

\begin{Proposition}\label {1.2}  If $N$ is
a $kG$-Hopf bimodule, then there exist a Hopf quiver $(Q, G,r )$, an
${\rm RSR} (G, r, \overrightarrow \rho, u)$ and a $kG$-Hopf bimodule
$(kQ_1^c, \alpha ^-, \alpha ^+)$ with \begin {eqnarray} \label
{e1.10.11} \alpha ^- (h \otimes a^{(i, j)}_{y,x}) := h\cdot a^{(i,
j)}_{y,x}=a^{(i, j)}_{hy,hx},\ \ \alpha ^+  ( a^{(i,
j)}_{y,x}\otimes h) := a^{(i, j)}_{y,x}\cdot h= \sum _{s\in J_C(i)}
k_{C, h}^{(i, j, s)}a^{(i, s)}_{yh,xh} \end {eqnarray} for some
$k_{C, h} ^{(i, j, s)}\in k$ such that $N \cong (kQ_1^c, \alpha ^-,
\alpha ^+)$ as $kG$-Hopf bimodules, where $x, y, h\in G$ with
$x^{-1}y=g^{-1}_{\theta}u(C)g_{\theta}$, $\zeta_{\theta}$ is given
by  (0.3) in \cite{ZZC04}, $C\in{\mathcal K}_r(G)$, $i\in I_C(r,
u)$, $j\in J_C(i)$, $x_C^{(i,j)} \cdot \zeta_{\theta}(h) = \sum_{s
\in J_C(i)} k_{C,h}^{(i, j, s)}x_C^{(i,s)}$.

\end {Proposition}
\noindent {\bf Proof.} Since $N$  is a $kG$-Hopf bimodule, there exists an
object \ \ \  \ \ \ \ $\prod_{C \in {\mathcal K}(G)} M(C) $ in
$\prod_{C \in {\mathcal K}(G)} {\mathcal M}_{kZ_{u(C)}}$ such that
$M(C)$ is a $kZ_{u(C)}$-module for  any $C \in {\mathcal K }(G)$ and
$N \cong \bigoplus _{y =
xg_\theta ^{-1} u(C) g_\theta, \  x, y \in G}\   x \otimes
M(C)\otimes _{kZ_{u(C)}} g_{\theta}$ as $kG$-Hopf bimodules by \cite
{CR97} or Theorem 1 in \cite{ZZC04}. Let $r = \sum _{C \in {\mathcal K}
(G)}r_C C$ with $ r_C ={\rm dim}M(C)$ for any $C \in {\mathcal K}
(G)$. Notice that ${\rm dim}M(C)$ denotes the  cardinal number of a
basis of $M(C)$ when $M(C)$ is infinite dimensional.  Since $M(C)$
is a $kZ_{u(C)}$-module and $k Z_{u(C)}$ is semisimple, there exists
a family of irreducible representations $\{(X_C^{(i)}, \rho
_C^{(i)}) \mid i \in I_C(r, u)\}$ such that $M(C)=\bigoplus_{i\in
I_C(r, u)}(X_C^{(i)},\rho_C^{(i)})$.
 Let $\{{x_C^{(i,j)}\mid j \in J_C(i)}\}$ be a k-basis of $X_C^{(i)}$ for
any $i\in I_C(r, u)$. Then for any $h\in G$ there are some
$k_{C,h}^{(i, j, s)} \in k$ such that $x_C^{(i,j)} \cdot
\zeta_{\theta}(h) = \sum_{s \in J_C(i)} k_{C,h}^{(i, j,
s)}x_C^{(i,s)}$ for all $i\in I_C(r, u)$ and $j \in J_C(i)$ since
$x_C^{(i,j)} \cdot \zeta_{\theta}(h) \in X_C^{(i)}$.

It remains to show that $(kQ_1^c, \alpha ^-, \alpha ^+)$ is
isomorphic to $\bigoplus _{y = xg_\theta ^{-1} u(C) g_\theta, \  x,
y \in G}\ x \otimes M(C)\otimes _{kZ_{u(C)}} g_{\theta}$ as
$kG$-Hopf bimodules. Observe that there is a canonical
$kG$-bicomodule isomorphism $\varphi: kQ_1\rightarrow \bigoplus _{y
=x g_\theta ^{-1} u(C) g_\theta, \  x, y \in G}\   x \otimes
M(C)\otimes _{kZ_{u(C)}} g_\theta$ given by
\begin{eqnarray}\label{indbya}
\varphi(a^{(i,j)}_{y,x})=x\otimes
x_C^{(i,j)}\otimes_{kZ_{u(C)}}g_{\theta}
\end{eqnarray}
 where $x, y\in G$ with $x^{-1}y=g^{-1}_{\theta}u(C)g_{\theta}$,
$C\in{\mathcal K}_r(G)$ , $i\in I_C(r, u)$ and $j\in J_C(r)$. We have
\begin{eqnarray*}
\varphi (\alpha ^- (h \otimes a^{(i,j)}_{y,x})) &=& \varphi (
a^{(i,j)}_{hy,hx})=
hx \otimes x^{(i,j)}_C \otimes_{kZ_{u(C)}}g_{\theta}  \ \ \\
&=&    h \cdot (x \otimes x^{(i,j)}_C \otimes_{kZ_{u(C)}}g_{\theta})
\ \ \ (\hbox {see }  (1.2) ~in ~\cite{ZZC04} ) \\
&=& h \cdot  \varphi ( a^{(i,j)}_{y,x}) \ \ .
\end{eqnarray*}
Thus  $\varphi$ is a left $kG$-module isomorphism.  Moreover
\begin{eqnarray*}
 \alpha ^+(\varphi (a^{(i,j)}_{y,x}) \otimes h) &=& xh
\otimes x_C^{(i,j)} \cdot \zeta_{\theta}(h) \otimes g_{\theta'}  \\
&=&xh \otimes(\sum_{s \in J_C(i)}{k_{C,h}^{(i, j, s)}}x_C^{(i,s)})
\otimes g_{\theta'}) \\
&=& \varphi (\sum_{s \in J_C(i)}{k_{C,h}^{(i,j, s)}}a_{yh, xh}^{(i,s)})\\
&=& \varphi ( \alpha ^+ (a_{y, x}^{(i,j)} \otimes h))  \ \ ( \hbox
{by } (\ref {e1.10.11})).
\end{eqnarray*}
Consequently, $\varphi$ is a $kG$-Hopf bimodule isomorphism. $\Box$

Let $(kQ_1^c, G, r, \overrightarrow \rho, u)$ denote the $kG$-Hopf
bimodule $(kQ_1^c, \alpha ^-, \alpha ^+)$ given in Proposition \ref
{1.2}. Furthermore, if $(kQ_1^c, kQ_1^a)$ is an arrow dual pairing, i.e. $kQ_1^c$ is isomorphic to
the dual of  $kQ_1^a$ as $kG$-Hopf bimodules or  $kQ_1^a$ is isomorphic to
the dual of  $ kQ_1^c$  as $(kG)^*$-Hopf bimodules under the isomorphisms
in Lemma 1.7 of \cite{ZZC04} (c.f. the argument before Definition 1.8 of \cite{ZZC04}),
then we denote the $(kG)^*$-Hopf bimodule  $kQ_1^a $ by $(kQ_1^a, G,
r, \overrightarrow \rho, u)$. We obtain six quiver Hopf algebras
$kQ^c (G, r, \overrightarrow \rho, u)$, $kQ^s (G, r, \overrightarrow
\rho, u)$, $kG[ kQ_1^c, G, r, \overrightarrow \rho, u]$, $kQ^a (G,
r, \overrightarrow \rho, u)$, $kQ^{sc} (G, r, \overrightarrow \rho,
u),$ $(kG)^*[ kQ_1^a, G, r, \overrightarrow \rho,$ $ u]$, called the
quiver Hopf algebras determined by ${\rm RSR} (G,$ $ r,$ $
\overrightarrow \rho,u)$.

From Proposition \ref{1.2}, it seems that the right $kG$-action on
$(kQ_1^c,G, r, \overrightarrow \rho, u)$ depends on the choice of
the set $\{g_{\theta}\mid \theta\in\Theta_C\}$ of coset
representatives of $Z_{u(C)}$ in $G$ (see, (\ref {e0.1})  or  (0.1) \cite{ZZC04}).
The following lemma shows that $(kQ_1^c,G, r, \overrightarrow \rho,
u)$ is, in fact, independent of the choice of the coset
representative set $\{g_{\theta}\mid \theta\in\Theta_C\}$, up to
$kG$- Hopf bimodule isomorphisms. For a while, we write $(kQ_1^c,G,
r, \overrightarrow \rho, u)=(kQ_1^c,G, r, \overrightarrow \rho, u,
\{g_{\theta}\})$ given before. Now let $\{h_{\theta}\in G\mid
\theta\in\Theta_C\}$ be another coset representative set of
$Z_{u(C)}$ in $G$ for any $C\in{\mathcal K}(G)$. That is,
\begin{eqnarray}\label{ncosetde}
G=\bigcup_{\theta\in\Theta_C}Z_{u(C)}h_{\theta}.
\end{eqnarray}
\begin{Lemma}\label{1.3}
With the above notations, $(kQ_1^c,G, r, \overrightarrow \rho, u,
\{g_{\theta}\})$ and $(kQ_1^c,G, r, \overrightarrow \rho, u, $
$\{h_{\theta}\})$ are isomorphic $kG$-Hopf bimodules.
\end{Lemma}

\noindent {\bf Proof.} We may assume \ \
$Z_{u(C)}h_{\theta}=Z_{u(C)}g_{\theta}$ \ \ for any \ $C\in{\mathcal
K}(G)$ and \ \  $\theta\in\Theta_C$. Then \ \ \
$g_{\theta}h^{-1}_{\theta}\in Z_{u(C)}$. Now let \ $x, y, h\in G$
with $x^{-1}y=g^{-1}_{\theta}u(C)g_{\theta}$. Then $x^{-1}y =$\ $
h^{-1}_{\theta}(g_{\theta}h^{-1}_{\theta})^{-1}u(C)(g_{\theta}h^{-1}_{\theta})h_{\theta}
$ \ $=h^{-1}_{\theta}u(C)h_{\theta}$ and
\begin{eqnarray}\label {e1.131}
h_{\theta}h=(h_{\theta}g^{-1}_{\theta})g_{\theta}h =
(h_{\theta}g^{-1}_{\theta})\zeta_{\theta}(h)g_{\theta'}=
(h_{\theta}g^{-1}_{\theta})\zeta_{\theta}(h)(g_{\theta'}h^{-1}_{\theta'})h_{\theta'}\ ,
\end {eqnarray} where $g_{\theta}h=\zeta_{\theta}(h)g_{\theta'}$.

For any $C \in {\mathcal K}(G)$, let $M(C)$ be a right
$kZ_{u(C)}$-module. Let  $N:= \bigoplus _{y = xg_\theta ^{-1} u(C)
g_\theta, \ x, y \in G}\ x \otimes M(C)\otimes _{kZ_{u(C)}}
g_{\theta}$ and $M:= \bigoplus _{y = xh_\theta ^{-1} u(C) h_\theta, \
x, y \in G}\ x \otimes M(C)\otimes _{kZ_{u(C)}} h_{\theta}$  be two
$kG$-Hopf bimodules. It is sufficient to show $N\cong M$ as
$kG$-Hopf bimodules by the proof of Proposition  \ref {1.2}.

Considering $x \otimes w \otimes _{k Z_{u(C)}} g_{\theta}  =x \otimes w \cdot g_\theta h_\theta ^{-1} \otimes _{k Z_{u(C)}}
h_\theta $, we have that $f: N\rightarrow M$ given by
$$ f ( x \otimes w \otimes _{k Z_{u(C)}} g_\theta ) =
x \otimes w \cdot g_\theta h_\theta ^{-1} \otimes _{k Z_{u(C)}}
h_\theta$$ for any $w\in M(C)$,  any $x, y\in G$ with
$x^{-1}y=g^{-1}_{\theta}u(C)g_{\theta}$, $C\in{\mathcal K}_r(G)$ and
$i\in I_C(r, u)$, is a $k$-linear isomorphism. It is clear that $f$ is  a $kG$-bicomodule isomorphism and a left
$kG$-module isomorphism from $N$ to $M$. Moreover
\begin{eqnarray*}
&&(f(x\otimes x_C^{(i,j)}\otimes_{kZ_{u(C)}}g_{\theta}))\cdot h
\\
&=& (x\otimes (x_C^{(i,j)}) \rho _C^{(i)} (g _\theta h_\theta ^{-1}
)\otimes_{kZ_{u(C)}}h_{\theta}) \cdot h\\
&=& xh\otimes (x_C^{(i,j)}) \rho _C^{(i)} ( \zeta _\theta (h)  )
\rho _C^{(i)} (g _{\theta '}h_{\theta '} ^{-1}
)\otimes_{kZ_{u(C)}}h_{\theta'}  \ \ (\hbox {by (\ref {e1.131})})\\
&=& f((x\otimes x_C^{(i,j)}\otimes_{kZ_{u(C)}}g_{\theta})\cdot h ),
\end {eqnarray*} for any  $x, y, h\in G$, $i\in I_C(r, u)$,  $j \in J_C(i)$,
$C\in{\mathcal K}_r(G)$  with
$x^{-1}y=g^{-1}_{\theta}u(C)g_{\theta}$. Thus $f$ is a right $kG$-module homomorphism.
\ \ $\Box$

Now we state one of  our main results, which classifies the quiver Hopf algebras.

\begin {Theorem} \label {1}
  Let $(G, r, \overrightarrow {\rho }, u)$ and
$(G', r', \overrightarrow{ \rho '}, u')$ be two {\rm RSR}s. Then
the following statements are equivalent:

{\rm(i)} ${\rm RSR} (G, r , \overrightarrow{\rho }, u)$ $\cong $
${\rm RSR} (G', r', \overrightarrow{\rho'}, u')$.

{\rm(ii)} There exists a Hopf algebra isomorphism $\phi:
kG\rightarrow kG'$ such that $(kQ_1^c, G, r , \overrightarrow{
\rho}, u) \cong\ _{\phi} ^{\phi ^{-1}}
 ( ( kQ_1' {}^c, G', r', \overrightarrow{\rho'}, u') ){}_{\phi} ^{\phi ^{-1}}$ as $kG$-Hopf bimodules.

{\rm(iii)} $kQ^c(G, r, \overrightarrow \rho, u)\cong k{Q'}^c(G, r,
\overrightarrow \rho, u)$.

 {\rm(iv)}  $kQ^s(G, r, \overrightarrow \rho, u) \cong
k{Q'}^s(G', r', \overrightarrow \rho', u')$.

 {\rm(v)} \
$kG[kQ_1^c,G, r, \overrightarrow \rho, u] $ $\cong kG'[kQ_1'{}^c,G', r',
 \overrightarrow \rho', u']$.\\
Furthermore, if $Q$ is finite, then the above are equivalent to the
following:

{\rm(vi)} $kQ^a(G, r, \overrightarrow \rho, u) $ $ \cong k{Q'}^a(G',
r', \overrightarrow \rho', u')$.

{\rm(vii)} $kQ^{sc}(G, r, \overrightarrow \rho, u) \cong
k{Q'}^{sc}(G', r', \overrightarrow \rho', u')$.

{\rm(viii)} \ $(kG)^*[kQ_1^a, G, r, \overrightarrow \rho, u]\cong
(kG')^*[kQ_1'{}^a, G', r', \overrightarrow \rho', u']$.
\end {Theorem}

\noindent {\bf Proof.}  By  Lemma 1.5 and Lemma 1.6 in \cite{ZZC04}, we only have
to prove {\rm (i)} $\Leftrightarrow$ {\rm (ii)}.

{\rm (i)} $\Rightarrow$ {\rm (ii)}. Assume that ${\rm RSR} (G, r ,
\overrightarrow{\rho}, u)$ $\cong $ ${\rm RSR} (G', r',
\overrightarrow{\rho'}, u')$. Let $(X_C^{(i)}, \rho _C^{(i)})$ and
$(X'{}_{C'}^{(i')}, \rho '{} _{C'}^{(i')})$ be irreducible
representations over $Z_{u(C)}$ and $Z_{u'(C')}$, respectively. Then
there exist a group isomorphism $\phi: G\rightarrow G'$,  an element
$h_C\in G$ such that $\phi(h^{-1}_Cu(C)h_C)=u'(\phi(C))$ for any
$C\in{\mathcal K}(G)$ and a bijective map $\phi_C: I_C(r,
u)\rightarrow I_{\phi (C)}(r', u')$ such that $(X^{(i)}_C,\rho_C^{(i)})
 \stackrel {\xi _C^{(i)}}{\cong} ( X'{}_{\phi (C)}^{(\phi _C(i))}, \rho'{}_{\phi(C)}^{(\phi_C(i))} \phi _{h_c})$ as right
$kZ_{u(C)}$-modules for all $i\in I_C(r, u)$.

 Now let
$G=\bigcup_{\theta\in\Theta_C}Z_{u(C)}g_{\theta}$ as in (\ref {e0.1})  or  (0.1) \cite
{ZZC04} for any $C\in{\mathcal K}(G)$. Let $x, y, h\in G$
with $x^{-1}y=g^{-1}_{\theta}u(C)g_{\theta}$, $C\in{\mathcal
K}_r(G)$ and $\theta\in\Theta_C$.  Assume that
$g_{\theta}h^{-1}_C=\zeta_{\theta}(h^{-1}_C)g_{\eta}$,
$g_{\theta}h=\zeta_{\theta}(h)g_{\theta'}$,
$g_{\eta}(h_Chh^{-1}_C)=\zeta_{\eta}(h_Chh^{-1}_C)g_{\eta'}$ and
$g_{\theta'}h^{-1}_C=\zeta_{\theta'}(h^{-1}_C)g_{\theta''}$ with
$\zeta_{\theta}(h^{-1}_C),$ $ \zeta_{\theta}(h), $ $
\zeta_{\eta}(h_Chh^{-1}_C), $ $\zeta_{\theta'}(h^{-1}_C)\in
Z_{u(C)}$ and $\eta, \theta', \eta', \theta''\in \Theta_C$. Then we have
\begin{eqnarray}\label {e1.106}
g_{\theta}hh^{-1}_C&=&\zeta_{\theta}(h) g_{\theta'}h^{-1}_C\ =
\zeta_{\theta}(h)\zeta_{\theta'}(h^{-1}_C)g_{\theta''}
\end{eqnarray} and
\begin{eqnarray}\label {e1.101}
g_{\theta}hh^{-1}_C&=&(g_{\theta}h^{-1}_C)(h_Chh^{-1}_C)=\zeta_{\theta}(h^{-1}_C)g_{\eta}(h_Chh^{-1}_C)\
=\zeta_{\theta}(h^{-1}_C)\zeta_{\eta}(h_Chh^{-1}_C)g_{\eta'}.
\end{eqnarray}
It follows that
\begin{eqnarray}\label{coeff}\theta''=\eta' \ \ \ \hbox {and } \ \
\ \zeta_{\theta}(h)\zeta_{\theta'}(h^{-1}_C)=
\zeta_{\theta}(h^{-1}_C)\zeta_{\eta}(h_Chh^{-1}_C).
\end{eqnarray}
 Moreover, we have
$(xh)^{-1}(yh)=h^{-1}g^{-1}_{\theta}u(C)g_{\theta}h=
g^{-1}_{\theta'}u(C)g_{\theta'}$ and
$g_{\theta}=g_{\theta}h^{-1}_Ch_C=\zeta_{\theta}(h^{-1}_C)g_{\eta}h_C$.
Thus $$\begin{array}{rcl}
\phi(x)^{-1}\phi(y)&=&\phi(x^{-1}y)=\phi(g^{-1}_{\theta}u(C)g_{\theta})
=\phi(h^{-1}_Cg^{-1}_{\eta}u(C)g_{\eta}h_C)\\
&=& \phi(h^{-1}_Cg^{-1}_{\eta}h_C)\phi(h^{-1}_Cu(C)h_C)\phi(h^{-1}_Cg_{\eta}h_C)\\
&=&\phi(h^{-1}_Cg_{\eta}h_C)^{-1}u'(\phi(C))\phi(h^{-1}_Cg_{\eta}h_C).
\end{array}$$
We also have
\begin{eqnarray}\label {e1.1311}
\phi(h^{-1}_Cg_{\eta}h_C)\phi(h)&=&\phi(h^{-1}_Cg_{\eta}h_Chh^{-1}_Ch_C) \nonumber \\
&=&\phi(h^{-1}_C\zeta_{\eta}(h_Chh^{-1}_C)g_{\eta'}h_C) \nonumber \\
&=&\phi(h^{-1}_C\zeta_{\eta}(h_Chh^{-1}_C)h_C)\phi(h^{-1}_Cg_{\eta'}h_C).
\end {eqnarray}
Since
$$\begin{array}{rcl} \phi _{h_C} (g_{\eta '})
 \phi (xh)^{-1} \phi (yh) &=&
 \phi (h_C^{-1}( \zeta _\eta (h_C hh_C^{-1}))^{-1} (\zeta _\theta (h_C^{-1}))^{-1}
 u(C)g_\theta h)  \ \ \ (  \hbox {by } \ref {e1.101}) \\
 &=& u'(\phi (C) )\phi_{h_C}(g_{\eta'}) \ \ \ (  \hbox {by } \ref {e1.101}), \\
\end{array}$$
\begin{eqnarray}\label {e1.102}
 \phi (xh)^{-1} \phi (yh)  &=& \phi _{h_C} (g_{\eta '}) ^{-1}u'(\phi (C))\phi_{h_C}(g_{\eta'}).
\end{eqnarray}

It is clear
\begin{eqnarray}\label{coset}
G'=\bigcup_{\theta\in\Theta_C}Z_{u'(\phi(C))}(\phi(h^{-1}_Cg_{\theta}h_C))
\end{eqnarray}
is a coset decomposition of $Z_{u'(\phi(C))}$ in $G'$ for any
$\phi(C)\in{\mathcal K}(G')$.

Let\begin{align*} N:&= \bigoplus _{y = xg_\theta ^{-1} u(C)
g_\theta, \ x, y \in G}\ x \otimes N(C)\otimes _{kZ_{u(C)}}
g_{\theta} \ \ \ \ \ \hbox { and } \\M:&= \bigoplus _{\phi (y) =
\phi (x) \phi _{h_C}(g_\eta ^{-1}) u'(\phi (C)\phi _{h_C}(g_\eta ),
\ x, y \in G}\ \phi (x) \otimes M(\phi (C))\otimes _{kZ_{u'(\phi
(C))}} \phi _{h_C}(g_\eta )\end{align*} with $N(C) := \bigoplus _{i
\in I_C(r, u) }X_C^{(i)}$ and $M(\phi (C)) := \bigoplus _{i\in
I_C(r, u)} X'{}_{\phi(C)}^{(\phi _C (i))}$. It suffices to show
$N\cong \ ^{\phi ^{-1}} _\phi M^{\phi ^{-1}}_\phi $ as $kG$-Hopf
bimodules by the proof of Proposition \ref {1.2}.

Considering  $x \otimes (w)\rho ^{(i)}_C ( (\zeta _\theta (h_C^{-1}))^{-1})\otimes _{kZ_{u(C)}}
g_{\theta} =  x \otimes w \otimes _{kZ_{u(C)}}
g_\eta h_c$, we have that $\psi: N \rightarrow M$ given by
$$\psi( x \otimes (w)\rho ^{(i)}_C ( (\zeta _\theta (h_C^{-1}))^{-1})\otimes _{kZ_{u(C)}}
g_{\theta})= \phi (x) \otimes \xi _C^{(i)}(w) \otimes _{kZ_{u'(\phi
(C))}} \phi _{h_C}(g_\eta )$$ for
 any $x, y\in G$ with $x^{-1}y=g^{-1}_{\theta}u(C)g_{\theta}$, and
$i\in I_C(r, u)$, $w \in X_C^{(i)}$,  where $C\in{\mathcal K}_r(G)$
and $g_{\theta}h^{-1}_C=\zeta_{\theta}(h^{-1}_C)g_{\eta}$ with
$\zeta_{\theta}(h^{-1}_C)\in Z_{u(C)}$ and $\theta, \eta\in
\Theta_C$, is a $k$-linear isomorphism. It is clear that $\psi$ is  a homomorphism not only of
$kG$-bicomodules from $N$ to $\ ^{\phi ^{-1}}  M^{\phi ^{-1}}$ but
also of left $kG$-modules from $N$ to $\  _\phi  M$.

For any $h \in G$ and $w \in X_C ^ {(i)}$, we have
\begin{align*} 
&\psi(( x \otimes (w)\rho ^{(i)}_C ( (\zeta _\theta (h_C^{-1}))^{-1}
)\otimes _{kZ_{u(C)}}  g_{\theta}     ) \cdot h)\\
 &=\psi(  ( x h\otimes (w)\rho ^{(i)}_C ( (\zeta _\theta
(h_C^{-1}))^{-1}  \zeta _\theta (h)  )\otimes _{kZ_{u(C)}}  g_{\theta'}    )\\
&=\phi (xh) \otimes \xi _C^{(i)}((w)\rho ^{(i)}_C ( (\zeta _\theta
(h_C^{-1}))^{-1}  \zeta _\theta (h) \zeta _{\theta '}(h_C^{-1}) ))
\otimes _{kZ_{u'(\phi (C))}} \phi _{h_C}(g_{\eta '} ) \ \ \ (\hbox { by (\ref {e1.102})})\\
&= \phi (x)\phi (h) \otimes ( \xi _C^{(i)}(w) )\rho ^{(\phi _C
(i))}_{\phi (C)} \phi _{h_C} ( (\zeta _\theta (h_C^{-1}))^{-1} \zeta
_\theta (h) \zeta _{\theta '} (h_C^{-1}) )) \otimes _{kZ_{u'(\phi (C))}}
\phi _{h_C}(g_{\eta '})\\
& \ \ (\hbox { by Definition \ref {1.1}})\\
 \end{align*}
 and
\begin{align*} &\psi
( x \otimes (w)\rho ^{(i)}_C ( (\zeta _\theta (h_C^{-1}))^{-1} )\otimes _{kZ_{u(C)}}
 g_{\theta}     ) \cdot \phi (h)\\
&=(\phi (x) \otimes \xi _C^{(i)}(w) \otimes _{kZ_{u'(\phi (C))}} \phi _{h_C}(g_\eta ) )\cdot \phi (h)\\
&= \phi (x) \phi (h)\otimes ( \xi _C^{(i)}(w) ) \rho ^{(\phi _C
(i))}_{\phi (C)} \phi _{h_C} ( \zeta _\eta (h_C h h_C^{-1}))\otimes _{kZ_{u'(\phi (C))}}
\phi _{h_C}(g_{\eta'} )   \ \ \ (\hbox { by (\ref {e1.1311})})    \\
&=  \phi (x) \phi (h)\otimes ( \xi _C^{(i)}(w) ) \rho ^{(\phi _C
(i))}_{\phi (C)} \phi _{h_C} ( (\zeta _\theta (h_C^{-1}))^{-1} \zeta
_\theta (h) \zeta _{\theta '}(h_C^{-1}) )) \otimes _{kZ_{u'(\phi
(C))}} \phi _{h_C}(g_{\eta'} ) \\
&  \ \ \ (\hbox { by (\ref {coeff})}),
 \end{align*}
which shows that $\psi$ is a right $kG$-module homomorphism.

{\rm (ii)} $\Rightarrow$ {\rm (i)}. Assume that there exist a Hopf algebra
isomorphism $\phi: kG\rightarrow kG'$ and a $kG$-Hopf bimodule
isomorphism $\psi: (kQ_1^c, G, r, \overrightarrow \rho,
u)\rightarrow\ ^{\phi^{-1}}_{\phi}(kQ_1'{}^c, G', r',
\overrightarrow \rho', u')^{\phi^{-1}}_{\phi}$. Then $\phi:
G\rightarrow G'$ is a group isomorphism. Let $C\in{\mathcal K}(G)$.
Then $\phi(u(C))$, $u'(\phi(C))\in\phi(C)\in{\mathcal K}(G')$, and hence
$u'(\phi(C))=\phi(h_C)^{-1}\phi(u(C))\phi(h_C)=\phi(h^{-1}_Cu(C)h_C)$
for some $h_C\in G$. Since $\psi$ is a $kG'$-bicomodule isomorphism
from $^{\phi}(kQ_1^c, G, r, \overrightarrow \rho, u)^{\phi}$ to
$(kQ_1'{}^c, G', r', \overrightarrow \rho', u')$ and
$\phi(h^{-1}_Cu(C)h_C)=u'(\phi(C))$, by restriction one gets a $k$-linear isomorphism
$$\psi_C:\ ^{h^{-1}_Cu(C)h_C}(kQ_1)^1\rightarrow\ ^{u'(\phi(C))}\! (kQ'_1)^1,\ x\mapsto \psi(x).$$
We also have a $k$-linear isomorphism
$$f_C:\ ^{u(C)}\! (kQ_1)^1\rightarrow\ ^{h^{-1}_Cu(C)h_C}(kQ_1)^1,\ x\mapsto h^{-1}_C\cdot x\cdot h_C.$$
Since $\phi(h^{-1}_Cu(C)h_C)=u'(\phi(C))$ and
$h^{-1}_CZ_{u(C)}h_C=Z_{h^{-1}_Cu(C)h_C}$, one gets
$\phi(h^{-1}_CZ_{u(C)}h_C)=Z_{u'(\phi(C))}$. Hence $\phi_{h_C}$ is
an algebra isomorphism from  $kZ_{u(C)}$ to $kZ_{u'(\phi(C))}$ 
sending $h$ to $\phi(h^{-1}_Chh_C)$. Using the hypothesis that
$\psi$ is a $kG$-bimodules homomorphism from  $(kQ_1^c, G, r,
\overrightarrow \rho, u)$   to $_{\phi}(kQ_1'{}^c,$\ $ G', r',
\overrightarrow \rho', u')_{\phi}$, one can easily check that the
composition $\psi_Cf_C$ is a right $kZ_{u(C)}$-module isomorphism
from $(^{u(C)}\! (kQ_1)^1, \lhd )$ to $((^{u'(\phi(C))}\! (kQ'_1)^1)_{\phi
_{h_C}}, \lhd )$. Indeed, for any $z \in ^{u(C)}\! (kQ_1)^1$, we have
\begin{align*}  \psi _Cf_C (z)\lhd \phi _{h_C} (h)  &= \psi _C (z \lhd h_C)\lhd \phi _{h_C}  (h) \\
& = \psi _C (z \lhd h_C\phi _{h_C} (h) ) \ \ \ (\hbox{since } \psi   \hbox { is a
bimodule homomorphism} )\\
&= \psi _Cf_C (z \lhd h).\\
 \end{align*}
Obviously, both  $^{u(C)}\! (kQ_1)^1$   and $(^{u'(\phi(C))}\!
(kQ'_1)^1)_{\phi _{h_C}}$ \ \ are semisimple right
$kZ_{u(C)}$-modules. Assume $^{u(C)}\! (kQ_1)^1 \cong \oplus _{i\in
I_C(r, u)} (X_C^{(i)}, \rho _C^{(i)})$ as right $kZ_{u(C)}$-modules
and $^{u'(\phi(C))}\! (kQ'_1)^1 \cong \oplus _{j\in I_{\phi
(C)}(r',u)} (X'{}_{\phi (C)}^{(j)}, \rho '{} _{\phi (C)}^{(j)})$ as
right $kZ_{u'(\phi (C))}$-modules, where $(X_C^{(i)}, \rho
_C^{(i)})$ is an irreducible right $kZ_{u(C)}$-module for any $i \in
I_C(r, u )$ and $(X'{}_{\phi (C)}^{(j)}, \rho '{} _{\phi (C)}^{(j)})$
is an irreducible  right $kZ_{u'(\phi (C))}$-module for any $j \in
I_{\phi (C)}(r', u')$. Therefore, there exists  a bijective map $\phi_C:
I_C(r, u)\rightarrow I_{\phi (C)}(r', u')$ such that $(X_C ^{(i)},
\rho_C^{(i)}) \cong (X'{}_{\phi (C)}^{(\phi _C(i))},
\rho_{\phi(C)}^{(\phi_C(i))} \phi _{h_C})$  as right
$kZ_{u(C)}-modules$ for all $i\in I_C(r, u)$. It follows that ${\rm
RSR} (G, r , \overrightarrow{\rho}, u)$ $ \cong $ $ {\rm RSR} (G',
r', \overrightarrow{\rho '},$ $ u').$ \ \ $\Box $

We have classified the  quiver Hopf algebras by means of
{\rm RSR}s. In other words,   ramification systems with irreducible
representations uniquely determine  their corresponding  quiver Hopf
algebras  up to graded Hopf algebra isomorphisms.

\begin{Proposition}\label {1.5} Let ${\rm RSR}(G, r, \overrightarrow \rho, u)$
and ${\rm RSR}(G, r, \overrightarrow {\rho'},
 u')$ be two {\rm RSR}s. If $u'(C) = h_C ^{-1}u(C) h_C$ and $\rho _C {}^{(i)} =
 \rho '{} _C ^{(i)} {\rm ad }^+_{h_C}
 $ with $I_C(r, u) = I_C(r, u')$ for any $C \in {\mathcal K}_r (G)$, $i \in I_C(r, u)$, where $h_C\in G$ and
 $ {\rm ad }^+ _{h_C} (g) = h_C^{-1} g h_C$,
 then  $ {\rm RSR}(G, r, \overrightarrow \rho, u)$ $\cong $ ${\rm RSR}(G, r, \overrightarrow {\rho'},
 u')$.
\end {Proposition}

\noindent {\bf Proof.} Let $\phi = id _G$ and $\phi _C = id_{I_C(r, u)}$ for
any $C\in {\mathcal K}_r (G).$  It is clear that $ {\rm RSR}(G, r,
\overrightarrow \rho, u)$ $\cong $ ${\rm RSR}(G, r, \overrightarrow
{\rho'}, u')$. $\Box$

 \noindent {\bf Remark.} This proposition means that the choice of map
 $u$ doesn't affect the classification of ${\rm RSR }$s. That is, if we fix a map
 $u_0$ from ${\mathcal K}(G)$ to $G$ with $u_0 (C) \in C$ for any $C \in {\mathcal
 K}(G)$, then for any $ {\rm RSR}(G, r,
\overrightarrow \rho, u)$, there exists $ {\rm RSR}(G, r,
\overrightarrow {\rho '}, u_0)$ such that $ {\rm RSR}(G, r,
\overrightarrow \rho, u)$ $\cong $ ${\rm RSR}(G, r, \overrightarrow {\rho'},  u_0)$.

\section{Classification of pointed Hopf algebras of type one }\label {s2}

A graded Hopf algebra $A = \oplus _{n=0} ^\infty A_{(n)}$ is said to be of  Nichols type, if
the diagram of $A$ is a Nichols algebra over $A_{(0)}$ (the definition
of diagram was given in \cite {AS98b} and Subsection 3.1 of \cite{ZZC04}).
Furthermore, if the coradical of $A$ is a group algebra, then $A$ is
called a pointed Hopf algebra of Nichols type.

For an ${\rm RSR}(G, r, \overrightarrow \rho, u)$ and a $kG$-Hopf
bimodule $(kQ_1^c, G, r, \overrightarrow {\rho}, u)$ with the module
operations $\alpha^-$ and $\alpha^+$, define a new left $kG$-action on $kQ_1$ by
$$g\rhd x:=g\cdot x\cdot g^{-1},\ g\in G, x\in kQ_1,$$
where $g\cdot x=\alpha^-(g\otimes x)$ and $x\cdot
g=\alpha^+(x\otimes g)$ for any $g\in G$ and $x\in kQ_1$. With this
left $kG$-action and the original left (arrow) $kG$-coaction
$\delta^-$, $kQ_1$ is a  Yetter-Drinfeld   $kG$-module. Let  $Q_1^1:=\{a\in Q_1 \mid
s(a)=1\}$, the set of all arrows with starting vertex $1$. It is
clear that $kQ_1^1$ is a  Yetter-Drinfeld   $kG$-submodule of $kQ_1$, denoted by
$(kQ_1^1, ad(G, r, \overrightarrow {\rho}, u))$.

\begin{Lemma}\label {2.1} {\rm (i)} If $H$ is a Hopf algebra with bijective antipode and
 $(B, \alpha _B^-, \delta
_B^-)$ is a graded braided Hopf algebra   in $^{H}_{H} {\mathcal YD}$
with $B_{(0)} =k1_B$, then ${\rm diag } (B\# H) = B \# 1_H \cong B$
as graded braided Hopf algebras in $^{H}_{H} {\mathcal YD}$.

{\rm (ii)}  $A$ is a pointed Hopf algebra of Nichols type if and
only if $A$ is isomorphic to the biproduct ${\mathfrak B}(V) \# kG $ as
graded Hopf algebras with the Nichols algebra ${\mathfrak B}(V)$
 over the group algebra  $kG$,  $A_{(0)} =kG$ and $ A_{(1)} = V \# kG$.
\end {Lemma}

\noindent {\bf Proof.} {\rm (i)} Obviously, $ {\rm diag} (B\# H) = B \otimes
1_H$. Define a map $\psi $ from  $ {\rm diag } (B\# H)$ to $B$ by
sending $x \otimes 1_H$ to $x$ for any $x\in B$. It is easy to check
that $\psi$ is a graded braided Hopf algebra isomorphism in
$^{H}_{H} {\mathcal YD}$.

{\rm (ii)} If $A$ is a pointed Hopf algebra of Nichols type, then
${\rm diag} (A)= {\mathfrak B}(V) $ is a Nichols algebra in
$^{kG}_{kG} {\mathcal YD}$ and the coradical  of $A$ is the group
algebra $kG$. Therefore $A \cong {\mathfrak B}(V) \# kG$ as graded
Hopf algebras. By Lemma 2.5 in \cite{AS98b}, $A_{(0)}= kG$ and
$A_{(1)} =  V \# kG.$

Conversely, clearly, $ {\rm diag } ({\mathfrak B} (V) \#kG)  = $
$({\mathfrak B} (V) \#kG)^{{\rm co }  kG} = $ $  {\mathfrak B} (V)\#
1 \cong $ $ {\mathfrak B} (V) $ by { \rm (i)} and the coradical of
${\mathfrak B} (V) \#kG$ is $({\mathfrak B} (V) \#kG)_{(0)}=kG$.
$\Box$

\begin{Lemma}\label {2.2} If $H = kG$ is a group algebra and $M$ is
an $H$-Hopf bimodule, then  the pointed  Hopf algebra $H[M]$ of  type
one is a Hopf algebra of Nichols  type. In particular, a
one-type-co-path Hopf algebra $kG [kQ_1^c, G, r,
\overrightarrow{\rho}, u]$ is a pointed Hopf algebra of Nichols type
and ${\rm diag} (kG [kQ_1^c, G, r, \overrightarrow{\rho}, u])$ $=
{\mathfrak B} (kQ_1^1, {\rm ad}(G, r, \overrightarrow{\rho}, u) )$.

\end {Lemma}
\noindent {\bf Proof.} By Proposition \ref {1.2} there exists an ${\rm RSR}(G,
r \overrightarrow{\rho}, u)$ such that $M \cong (kQ_1^c, G, r
\overrightarrow{\rho}, u)$ as $kG$-Hopf bimodules. Thus  $kG [M]
\cong kG [kQ_1^c]$ as graded Hopf algebras by Lemma 1.6 in \cite{ZZC04}. 
Therefore, it is enough to show that $kG [kQ_1^c]$ is
of Nichols type. Let $A := kG [kQ_1^c]$ and $R:= {\rm diag} (kG
[kQ_1^c])$. Obviously, $R_{(0)} =k.$ Now we show $R_{(1)} = k Q_1^1$. Obviously,
$ k Q_1^1 \subseteq
R_{(1)}$. Let $\alpha = \sum _{p=1}^n k_pb^{(p)} \in R_{(1)}$,  where $b^{(p)}$ is an arrow from
$x^{(p)}$ to $y^{(p)}$ with $0 \not= k_p \in k$,  and $b^{(1)}, b^{(2)}, \cdots, b^{(n)} $  are different from each other.
Therefore $\sum _{p=1}^n k_pb^{(p)} \otimes 1 = \sum _{p=1}^n k_pb^{(p)} \otimes x^{p}$, which implies
$x^{(p)}=1 $ for $1\le p \le n.$ Thus $\alpha \in kQ_1^1.$ We next show
that $R$ is generated by $R_{(1)}$ as  algebras. Let $\mu$ denote the multiplication
and let $B$ denote the algebra generated by $kQ_1^1$ in $kG [kQ_1^c]$. Obviously,
$B \subseteq R.$ It follows from the argument in Subsection 3.1 of \cite{ZZC04} that
$\alpha^+ := \mu (id \otimes \iota _0)$ is an algebraic isomorphism from $R \# kG$ to $kG [kQ_1^c]$. For any $x, y \in G$
and any arrow $a_{y, x}$ from $x$ to $y$, we have
\begin {eqnarray*}
a_{y, x} & =& x \cdot a_{x^{-1}y, 1} = \mu (x \otimes a_{x^{-1}y, 1} ) = \mu ( \alpha ^+
(1 \# x) \otimes \alpha ^+(a_{x^{-1}y,1} \# 1))\\
 &=& \alpha ^+ (  \mu ( (1\# x ) \otimes (a_{x^{-1}y, 1}\#1))) =
\alpha^+ (x\rhd a_{x^{-1}y, 1}  \#x) \in \alpha^+ (B\#kG).
\end {eqnarray*} Therefore
$\alpha ^+ (B\# kG) = \alpha ^+ (R\# kG)$ and so $B= R.$

It is sufficient to show $P(R) = kQ_1^1$,  where $P(R)$ denotes the set of all
primitive elements in $R$. For any $a \in Q_1^1$ with $\delta ^- (a)
=y$ and $\delta ^+ (a)= 1$, we have $\Delta _R (a)=$ $(\omega \otimes
id) \Delta _A (a)=$ $ 1\otimes a +a \otimes 1$ (see Section 3 of \cite{ZZC04}), 
i.e. $ kQ_1^1 \subseteq P(R)$, where $\omega = \mu _A (id \otimes \iota _0\pi _0S)\Delta _A$.

Conversely, we shall show $P(R) \subseteq  kQ_1^1$ by the following
two steps. Obviously, $kG \cap P(R) =0$ and $P(R)$ is a graded subspace of $R$.

(i) Assume that  $\alpha = a_{x_nx_{n-1}} a_{x_{n-1}x_{n-2}} \cdots
a_{x_1x_0}$ is a path from vertex $x_0$, via arrows $a_{x_1x_0}$ $,
\cdots,$ $a_{x_{n-1}x_{n-2}}$, $a_{x_nx_{n-1}}$, to vertex $x_{n}$.
Then $\omega (\alpha ) = \alpha \cdot x_0^{-1}$.

{\rm (ii)} Let $v=\sum _{p =1}^m k_p\alpha _p \in P(R)$, where $\alpha _p
=b_{x_n x_{n-1} }^{(p)}b_{x_{n-1} x_{n-2} }^{(p)}\cdots b_{x_1 x_{0}
}^{(p)}$ is a path with $n>1$, $k_p \in k$ for $p =1, 2, \cdots, m$,
and $b_{x_j x_{j-1} }^{(p)}$ is an arrow from vertex $x_{j-1}$ to
vertex $x_j$ for $j=1,2,\cdots,n$. We shall show that $k_p =0$ for
$p =1, 2, \cdots, m.$ Indeed,
\begin {eqnarray*}
\Delta _R (v) &=&\sum _{p =1}^m \sum _{j=0}^n  k_p(b_{x_n x_{n-1}
}^{(p)}b_{x_{n-1} x_{n-} }^{(p)}\cdots b_{x_{j+1} x_{j} }^{(p)}
)\cdot (x^{(p)}_j)^{-1} \otimes b_{x_j x_{j-1}
}^{(p)}b_{x_{j-1} x_{j-2} }^{(p)}\cdots b_{x_1 x_{0} }^{(p)}\\
&=&\sum_{p =1}^mk_p(\alpha _p \otimes 1 + 1 \otimes \alpha _p).
\end {eqnarray*}
This implies
\begin {eqnarray} \label{e2.2.1}
\sum _{p =1}^m k_p(b_{x_n x_{n-1} }^{(p)}b_{x_{n-1} x_{n-2}
}^{(p)}\cdots b_{x_{j+1} x_{j} }^{(p)} )\cdot (x^{(p)}_j)^{-1}
\otimes b_{x_j x_{j-1} }^{(p)}b_{x_{j-1} x_{j-2} }^{(p)}\cdots
b_{x_1 x_{0} }^{(p)} &=& 0,
\end {eqnarray} for $j =1, 2, \cdots, n-1$, because of their length. For any $j$ with $1\leq
j\leq n-1$, assume that $  \{  p \mid  b_{x_j x_{j-1}
}^{(p)}b_{x_{j-1} x_{j-2} }^{(p)}\cdots b_{x_1 x_{0} }^{(p)} $ $ =
b_{x_j x_{j-1} }^{(1)}b_{x_{j-1} x_{j-2} }^{(1)}\cdots b_{x_1 x_{0}
}^{(1)} \}$ $ = \{1, 2, \cdots, m_1\}$ without loss of generality.
Therefore, by (\ref {e2.2.1}),
$$\sum _{p =1}^{m_1} k_p(b_{x_n x_{n-1} }^{(p)}b_{x_{n-1} x_{n-2} }^{(p)}\cdots b_{x_{j+1} x_{j}
}^{(p)} )\cdot (x^{(p)}_j)^{-1} \otimes b_{x_j x_{j-1}
}^{(p)}b_{x_{j-1} x_{j-2} }^{(p)}\cdots b_{x_1 x_{0} }^{(p)} =0,$$ since $ b_{x_j x_{j-1}
}^{(q)}b_{x_{j-1} x_{j-2} }^{(q)}\cdots b_{x_1 x_{0} }^{(q)} \not=
b_{x_j x_{j-1} }^{(1)}b_{x_{j-1} x_{j-2} }^{(1)}\cdots b_{x_1 x_{0}
}^{(1)}$  for   $q =  m_1+1, \cdots, m.$
This implies $\sum _{p =1}^{m_1} k_p b_{x_n x_{n-1}
}^{(p)}b_{x_{n-1} x_{n-2} }^{(p)}\cdots b_{x_{j+1} x_{j} }^{(p)} =0$
and $k_p =0$ for $p =1, 2, \cdots, m_1$. Similarly, we can show $k_p
=0$ for $p =m_1+1, \cdots, m.$ $\Box$

By Theorem 4.3.2 in \cite{BD}, the category $^{H}_{H} {\mathcal YD}$ of
Yetter-Drinfeld modules is equivalent to the category $^{H}_{H}
\!{\mathcal M}^{H}_{H}$ of $H$-Hopf bimodules, where $H$ is a Hopf
algebra with bijective antipode. Let $T$ and $U$ be the two
corresponding functors. For any
 $N \in ^{H}_{H} \!{\mathcal YD}$, according to Proposition 4.2.1 in \cite{BD}, $T(N) :=
 N\rtimes H=N\otimes H$ as vector spaces, and the
actions and coactions are given as follows: the left (co)actions are
diagonal and right (co)actions are induced by $H$. Explicitly, $g
\cdot (x \otimes h) := g \cdot x \otimes gh$, \ $(x \otimes h) \cdot
g = x\otimes hg$, \ $\delta^-_{N\rtimes H} (x\otimes h):= \sum _{x}
x_{(-1)} h \otimes x_{(0)} \otimes h; $  \ $\delta ^+ _{N\rtimes H}
(x\otimes h):= x \otimes h \otimes h$, where $\delta _N ^-(x) = \sum
_{x} x_{(-1)} \otimes x_{(0)}$, $x\in N, h, g \in H.$ For any $M \in
^{H}_{H} \! {\mathcal M}^{H}_{H}$, according to Equations (7) and
(21) in \cite{BD}, $U(M)$ is the coinvariant of $M$ as a vector space, i.e.,
$U(M) : = M ^{{\rm co } H} := \{ x \in M \mid \delta _N^+ (x) = x
\otimes 1\}$. The left action is left adjoint action and the left
coaction is the restricted coaction of the original coaction of $M$.
That is,
$$ \alpha _{U(M)}^- (h \otimes x ) = h \triangleright _{\rm ad }x :
=\alpha _M^+( \alpha _M^{-} (h \otimes x ) \otimes h^{-1}) =(h \cdot
x )\cdot h^{-1}$$ and  $ \delta _{U(M)} ^- (x) = \delta _M^-(x)$ for
any $h \in H, x\in U (M).$ In fact, $TU(M)=U(M)\rtimes H$ and $UT(N)
= N \otimes 1_H$. Let $\lambda_N$ be map from $N\otimes 1_H$ to
$N$ sending $x\otimes 1_H$ to $x$ for any $x\in N$, and let
$\nu_M$ be map from $U(M)\rtimes H$ to $M$ sending $x\otimes h$
to $\alpha _M^+(x \otimes h)=x\cdot h$ for any $x\in U(M)$ and $h\in
H$. $\lambda$ and $\nu$ are the natural isomorphisms from functor
$UT$ to $id$ and from functor $TU$ to $id$, respectively. Note that
the inverse of $\nu_M$ is $(\alpha^+_M \otimes id) (id \otimes
S\otimes id) (\delta ^+_M\otimes id) \delta_M^+.$

\noindent {\bf Remark.} We have $U  (kQ_1^c, G, r,
\overrightarrow{\rho}, u) = (kQ_1^1, {\rm ad }( G, r,
\overrightarrow{\rho}, u))$ by the proof of Lemma \ref {2.2}.
\begin{Lemma}\label{2.3}
Assume that $\phi$ is a Hopf algebra isomorphism from $H$ to $H'$.
Let $N\in{}^{H'}_{H'}{\mathcal YD}$ and $M\in{}^{H'}_{H'}\!{\mathcal
M}^{H'}_{H'}$. Then
$$ T ( {}{^\phi }^{-1}_\phi \! N) \cong  {}{ ^\phi} ^{-1}_\phi T(N)
{^\phi}^{-1} _\phi  \hbox {\ \ in \ } {}^{H}_{H} \!{\mathcal M}^{H}_{H}
\hbox { \ \ and \ \ }  U( {} {^\phi } ^{-1}_\phi M {^\phi }^{-1}
_\phi ) \cong {^ \phi} ^{-1}_\phi U(M) \hbox {\ \ in }\  {} ^{H}_{H}
{\mathcal YD}.$$
\end {Lemma}

\noindent {\bf Proof.} The first isomorphism is given by sending $x\otimes h$ to $x \otimes \phi (h)$
for any $x\in N, h\in H$; the second one is identity.  $\Box$

\begin{Proposition}\label {2.4}
{\rm (i)} If $N$ is a Yetter-Drinfeld $kG$-module, then there exists
an $ {\rm RSR}(G, r, $ $\overrightarrow {\rho}, u)$ such that $N
\cong (kQ_1^1, ad (G, r, \overrightarrow \rho, u))$ as
Yetter-Drinfeld $kG$-modules.

{\rm (ii)} If ${\mathfrak B}(N)$ is a Nichols algebra in
$^{kG}_{kG}{\mathcal YD}$, then there exists an $ {\rm RSR}(G, r,
\overrightarrow {\rho}, u)$ such that ${\mathfrak B}(N) \cong
{\mathfrak B}(kQ_1^1, ad (G, r, \overrightarrow \rho, u))$ as
graded braided Hopf algebras in $^{kG}_{kG}{\mathcal YD}$.
\end{Proposition}

\noindent {\bf Proof.} {\rm (i)} Since $T(N)$ is a $kG$-Hopf bimodule, it
follows from Proposition \ref{1.2} that there exists an ${\rm
RSR}(G, r, \overrightarrow {\rho}, u)$ such that $T(N) \cong
(kQ_1^c, G, r, \overrightarrow \rho, u)$ as $kG$-Hopf bimodules.
Thus, $N \cong UT(N) $ $ \cong U(kQ_1^c, G, r, \overrightarrow \rho,
u)= (kQ_1^1, ad (G, r, \overrightarrow \rho, u))$ as Yetter-Drinfeld
$kG$-modules by Equations (7) and (21) in \cite{BD}.

{\rm (ii)} This follows from {\rm (i)} and Corollary 2.3 in \cite{AS02}.
$\Box$

\begin{Lemma}\label{2.5} Assume that $\phi$ is a Hopf algebra isomorphism
from $H$ to $H'$. Let $R$ and $R'$ be graded braided Hopf algebras
in $^{H}_{H}{\mathcal YD}$ and $^{H'}_{H'}{\mathcal YD}$ with $R_{(0)}=k1_R$
and $R_{(0)}'=k1_{R'}$, respectively. If $R$ and $^{\phi^{-1}}_\phi
R'$ are isomorphic as graded braided Hopf algebras in $ {} ^{H}_{H}
{\mathcal YD}$, then biproducts $R \# H \cong R'\# H'$ as graded Hopf algebras.
\end{Lemma}

\noindent {\bf Proof.} Let $\psi$ be a graded braided Hopf algebra isomorphism
from $R$ to $^{\phi ^{-1}}_\phi R'$ in $^{H}_{H}{\mathcal YD}$. Define a
map $\nu$ from $R \#H$ to $R'\#H'$ by sending $x\otimes h$ to $\psi
(x)\otimes \phi (h)$ for any $x\in R, h\in H$. It is easy to check
that $\nu$ is an isomorphism of graded Hopf algebras. $\Box$

\begin{Theorem}\label{2}
If $A$ is a pointed Hopf algebra of Nichols type with coradical
$kG$, a group algebra, then there exists a unique ${\rm RSR}(G, r ,
\overrightarrow{\rho }, u)$, up to isomorphism, such that $A\cong
kG[kQ_1^c, G, r, \overrightarrow{\rho}, u]$ as graded Hopf algebras.
\end{Theorem}

\noindent {\bf Proof.} By Lemma \ref {2.1}, $A \cong {\mathfrak B}(V) \# kG$
as graded Hopf algebras. By Proposition \ref {2.4} {\rm (ii)}, there
exists a ${\rm RSR} (G, r , \overrightarrow{\rho }, u)$ such that
${\mathfrak B}(V) \cong {\mathfrak B}(kQ_1^1, {\rm ad}(G, r ,
\overrightarrow{\rho }, u))$ as graded braided Hopf algebras in
$^{kG}_{kG} \! {\mathcal YD}$. Thus
\begin {eqnarray*}
kG [kQ_1^c, G, r , \overrightarrow{\rho }, u ] &\cong& {\mathfrak
B}(kQ_1^1,
{\rm ad}(G, r , \overrightarrow{\rho }, u)) \# kG \ \ (\hbox {by Lemma \ref {2.2}} )\\
&\cong &{\mathfrak B}(V)\# kG \ \ (\hbox {by Lemma \ref {2.5}})\\
&\cong & A.
\end{eqnarray*}
The uniqueness follows from Theorem 3 in \cite{ZZC04}. $\Box$

Considering  Theorem \ref {2} and Lemma \ref {2.2}, we have that $A$ is a
pointed Hopf algebra of Nichols type
if and only if $A$ is a pointed Hopf algebra of  type one.

\section {\bf Classification of Nichols Algebras }\label{s3}
\begin{Theorem}\label {3}
 Let $(G, r, \overrightarrow {\rho }, u)$ and
$(G', r', \overrightarrow{ \rho '}, u')$ be two {\rm RSR}s. Then
the following statements are equivalent:

{\rm(i)} ${\rm RSR} (G, r , \overrightarrow{\rho }, u)$ $\cong $
${\rm RSR} (G', r', \overrightarrow{\rho'}, u')$.

{\rm(ii)} There exists a Hopf algebra isomorphism $\phi:
kG\rightarrow kG'$ such that $(kQ_1^1, ad (G, r , \overrightarrow{
\rho}, u)) $ $  \cong\ $ $  _{\phi} ^{\phi ^{-1}}
  ( kQ_1' {}^1, ad ( G', r', \overrightarrow{\rho'}, u')) $ as Yetter-Drinfeld $kG$-modules.

{\rm(iii)} There is a Hopf algebra isomorphism $\phi: kG\rightarrow
kG'$ such that ${\mathfrak B}(kQ_1^1, {\rm ad}  (G, r ,
\overrightarrow{ \rho}, u)) $ $ \cong\ $ $_{\phi} ^{\phi ^{-1}}
{\mathfrak B}( kQ_1' {}^1, {\rm ad} ( G', r',
\overrightarrow{\rho'}, u'))$ as graded braided Hopf algebra in
$^{kG}_{kG} {\mathcal YD}$.
\end{Theorem}

\noindent {\bf Proof.} ${\rm (i)} \Rightarrow {\rm (ii)}$. We have
 \begin{eqnarray*}
  _{\phi} ^{\phi ^{-1}}
  ( kQ_1' {}^1, ad ( G', r', \overrightarrow{\rho'}, u')) &  = &  _{\phi} ^{\phi ^{-1}}
  U( kQ_1' {}^c, G', r', \overrightarrow{\rho'}, u')  \ \ \ (\hbox { by the remark before Lemma \ref {2.3}} )\\
&  \cong &  U ( _{\phi} ^{\phi ^{-1}}   (kQ_1' {}^c, G', r', \overrightarrow{\rho'}, u')_{\phi} ^{\phi ^{-1}}
)  \ \ \ (\hbox { by Lemma \ref {2.3}})\\
&\cong & U((kQ_1^c, G, r, \overrightarrow{\rho}, u)  ) \ \ \ (\hbox { by Theorem \ref {1}})\\
&=& (kQ_1^1, ad ( G, r, \overrightarrow{\rho}, u)).
\end{eqnarray*}

${\rm (ii)} \Rightarrow {\rm (iii)}$. By Lemma 2.7 in \cite{ZZC04},
${\mathfrak B} (kQ_1^1, {\rm ad}  (G, r , \overrightarrow{ \rho},
u)) $ $ \cong $ $ {\mathfrak B} ( ^{\phi ^{-1}} _\phi  ({kQ_1'}^1,
{\rm ad} ( G', r', \overrightarrow{\rho'}, u')))$ $ \cong $ $ \
^{\phi ^{-1}} _\phi \! {\mathfrak B} ( {kQ_1'}^1, {\rm ad} ( G', r',
\overrightarrow{\rho'}, u')) $ as graded braided Hopf algebras in
$^{kG}_{kG} {\mathcal YD}$.

${\rm (iii)} \Rightarrow {\rm (i)}$. ${\mathfrak B}(kQ_1^1) \# kG
\cong {\mathfrak B}(kQ_1'{}^1) \# kG'$ as graded Hopf algebras by
Lemma \ref {2.5}. Thus $kG [kQ_1^c] \cong {\mathfrak B}(kQ_1^1)\# kG
\cong {\mathfrak B} (kQ_1'{}^1)\#kG'\cong kG' [kQ'_1{}^c]$ as graded
Hopf algebras by Lemma \ref {2.2}. Now {\rm (i)} follows from
Theorem \ref {1}.
$\Box$\\

We have classified all   Nichols algebras by means of {\rm
RSR}s. In other words,   ramification systems with irreducible
representations uniquely determine  their corresponding Nichols
algebras up to pull-push graded braided  Hopf algebra isomorphisms.

\section {Classification of {\rm RSRs} over symmetric groups}

Let ${\rm ad}_{h}^-$ and ${\rm ad}_{h}^+$ denote the left and right
adjoint actions, respectively. That is, ${\rm ad}_{h}^-(x):= h x
h^{-1}$ for any $ x\in G.$ Let $\gamma _C = \mid \!Z_s \! \mid$ when
$s\in C\in {\mathcal K}(G)$. {\rm Aut}$G$ and {\rm Inn}$G$ denote
the automorphism group and inner automorphism group of $G$, respectively.

\begin{Definition}\label{5.1} Let
${\rm RSR} (G,   r,   \overrightarrow{\rho },   u)$ and ${\rm RSR} (G',   r',   \overrightarrow{\rho '},   u')$
be two RSRs. If there exists a bijective map $\phi _C$ from $I_C(r, u)$ to $I_C(r',u')$
such that $ \rho ^{(i)}_C \cong \rho '{}^{(\phi _C(i))}_C$ for any $i\in I_C(r, u)$
and $C\in {\mathcal K}_r(G)$ with  $G=G'$,   $r=r'$ and $u=u'$,  then
${\rm RSR} (G,   r,   \overrightarrow{\rho },   u)$ and ${\rm RSR} (G',   r',   \overrightarrow{\rho '},   u')$
are said to be of  the same type. Furthermore, if let  $\widehat {Z_{u(C)}}=\{ \xi _{u(C)}^{(i)} \mid i = 1,  2,
\cdots, \gamma _C\}$  and \ $n_C^{(i)} : =$ $ \mid \{ j \mid \rho _C
^{(j)} \cong \xi _{u(C)}^{(i)} \} \! \mid  $ \ for any $C\in
\mathcal {K}_r(G)$ \ \ and \ \ \ $1\le i \le \gamma _{u(C)}\  $,
then $\{(n_C^{(1)}, n_C^{(2)}, \cdots,  n_C^{(\gamma
_C)})\}_{C\in\mathcal {K}_r(G)}$ is called the type of ${\rm RSR}
(G,  r, \overrightarrow{\rho },  $ $ u)$.

\end{Definition}

\begin {Lemma} \label {5.5} If $\mathrm{Aut}G=\mathrm{Inn}G$, for example, $G= \mathbb S_n$ with
$n\not=6$, then  ${\rm RSR} (G,   r,   \overrightarrow {\rho},   u)$ and  ${\rm RSR} (G,   r,
\overrightarrow {\rho'},   u)$ are isomorphic if and only if they have the same type.

\end {Lemma}
\noindent {\bf Proof.} By 
Proposition 1.1 {\rm (ii)} in \cite{ZWW08}, 
$\mathrm{Aut}G=\mathrm{Inn}G$  when $G= \mathbb S_n$ with $n\not= 6$.

Let ${\rm RSR} (G',   r',   \overrightarrow {\rho'}, u')$ denote
${\rm RSR} (G,   r,   \overrightarrow {\rho'},   u)$ with $G= G'$,
$r=r'$ and $u=u'$ for convenience. If ${\rm RSR} (G, r,
\overrightarrow {\rho}, u)$ and ${\rm RSR} (G, r, \overrightarrow
{\rho'},   u)$ have the same type, then there exists a bijective map $\phi _C$ from $I_C(r, u)$ to
$I_{\phi (C)}(r', u')$ such that  $\rho'{} ^{(\phi _C(i))}_C\cong \rho
^{(i)}_C$ for any  $i \in I_C(r, u)$,   $C \in {\mathcal K}_r (G)$.
 Therefore they are isomorphic.

Conversely,  if  ${\rm RSR} (G,   r,   \overrightarrow {\rho},   u)$
and ${\rm RSR} (G,  r,  \overrightarrow {\rho'},   u)$ are
isomorphic, then there exist a $\phi \in {\rm Aut }(G)$, $h_C\in G$
and a  bijective map $\phi_C : I_C(r, u) \rightarrow I_{\phi
(C)}(r', u')$  such that $u'(\phi (C)) = \phi {\rm ad}^+_{h_C}
(u(C))$ and $\rho'{}_{\phi (C)}^{(\phi_C(i))}\phi {\rm ad}^+_ {h_C}
\cong \rho_C^{(i)}$ for any
 $C \in {\mathcal K}_r(G)$ and $i\in I_C(r, u)$. Since $\phi {\rm ad} _{h_C}^+ \in {\rm Aut } G $
 $= {\rm Inn} G$, there exists a $g_C \in G$ such that $\phi {\rm ad}
 _{h_C}^+ = {\rm ad} _{g_C}^+$. Therefore  $u(C)={\rm ad}^+_{g_{C}}(u(C))$ and
$\rho'{}_{C}^{(\phi_C(i))}{\rm ad}^+_ {g_C} \cong \rho_C^{(i)}$. That is,
  $g_C \in  Z_{u(C)}$ and  $\chi '{}_{C}^{(\phi_C(i))}{\rm ad}^+_{g_{C}}(h)$ $=\chi '{}_{C}^{(\phi_C(i))}(g_C ^{-1}h g_C
)$ $=\chi'{}_{C}^{(\phi_C(i))}(h)=\chi_C^{(i)}(h) $ for any $h \in
Z_{u(C)}$, where $\chi'{} _C^{(\phi _C(i))}$ and $\chi _C^{(i)}$
denote the characters of  $\rho'{}_{C}^{(\phi_C(i))}$ and $
\rho_C^{(i)}$, respectively. Consequently, $\chi'{}_{C}^{(\phi_C(i))}=\chi_C^{(i)}$ and
$\rho'{}_{C}^{(\phi_C(i))} \cong \rho_C^{(i)}$.  This implies that
${\rm RSR} (G, r, \overrightarrow {\rho},  u)$ and ${\rm RSR} (G, r,
\overrightarrow {\rho'},  u)$ have the same type. $\Box$

For a given ramification $r$ of $G$,  let  $ \Omega (G,  r)$ be the set of all
RSRs of G with the ramification $r$,  namely,  $ \Omega (G,  r) :
= \{ (G,  r,  \overrightarrow{ \rho},  u) \mid (G,   r,
\overrightarrow{\rho},  u) \mbox {\ is\ an  } {\rm RSR}\}$. Let ${\mathcal
N}(G,  r)$ be the number of isomorphism classes in $ \Omega (G,r)$.

\begin{Theorem}\label{5.6} Given a group $G$ and a ramification $r$ of $G$. 
Assume $\mathrm{Aut}G=\mathrm{Inn}G$, for example, $G= \mathbb S_n$
with $n\not=6$. Let  $u_0$ be a fixed map from ${\mathcal K}(G)
\rightarrow G$ with  $u_0(C) \in C$  for any $ C\in {\mathcal
K}(G)$. Let $\bar {\Omega }(G,   r,   u_0)$ denote the set
consisting of all elements with  distinct type in  $\{ (G,  r,
\overrightarrow{\rho }, u_0) \mid (G,  r,  \overrightarrow{\rho },
u_0)$ is an  {\rm RSR}  $\}$. Then $\bar {\Omega }(G,
 r,   u_0)$ becomes the representative system  of isomorphic classes in $\Omega (G,   r). $
\end{Theorem}
\noindent {\bf Proof.} For any ${\rm RSR}(G,  r,  \overrightarrow{\rho }, u)$,
there exists ${\rm RSR}(G,  r,  \overrightarrow{\rho '}, u_0)$ such
that ${\rm RSR}(G,  r,  \overrightarrow{\rho }, u)$ $ \cong {\rm
RSR}(G, r, \overrightarrow{\rho' }, u_0)$ by Proposition \ref {1.5}.
Using Lemma \ref {5.5}, we complete the proof. $\Box$

\begin {Corollary} \label {5.7}
 \begin{eqnarray*}
{\mathcal N} (G,  r)&=&\prod\limits_{C\in \mathcal {K}_r(G)} \tau _C,
\end {eqnarray*}
where $\tau _C = $ the number of elements in the set $ \{ (n_1, n_2, \cdots, n_{\gamma_C})
\in {\mathbb N} ^{\gamma_C} \mid
n_1  {\rm deg } \xi _1+ n_2{\rm deg } \xi _2 + \cdots +  n_{\gamma_C}{\rm deg } \xi _{\gamma_C}
= r_C\}$ for any $C\in {\mathcal K}_r (G)$.
\end {Corollary}
\noindent {\bf Remark.} For a given finite Hopf quiver $(Q, G, r)$ over the symmetric group $G={\mathbb S}_n$
with $n \not= 6$, every path algebra $ T_{(kG)^*} (kQ_1^a)$ over the quiver $(Q, G, r)$ admits exactly 
 ${\mathcal N} (G, r)$ non-isomorphic
graded Hopf algebra structures;  every path coalgebra $ T_{kG}^c (kQ_1^c)$ over the quiver $(Q, G,  r)$ admits
exactly  ${\mathcal N} (G, r)$ non-isomorphic graded Hopf algebra structures.

\section {Appendix}\label {S4}
We now consider the dual case of Theorem \ref {3}. If  $Q$ is finite, then $(kQ_1^a, G, r, \overrightarrow
{\rho}, u)$ is a $(kG)^*$-Hopf bimodule with comodule operations
$\delta^-$ and $\delta^+$. Define a new left $(kG)^*$-coaction on $kQ_1^a$ given by
$$\delta^-_{\rm coad}(x):=\sum_{x}x_{(-1)}S(x_{(0)(1)})\otimes x_{(0)(0)},\ \ \ \hbox { for any } x\in kQ_1^a,$$
i.e. adjoint coaction. With this left $(kG)^*$-coaction and the original left (arrow)
$(kG)^*$-action $\alpha^-$, $kQ_1^a$ is a  Yetter-Drinfeld   $(kG)^*$-module. Let
$kQ^{1a}_1$ denote the subspace spanned by $Q_1^1$ in $kQ_1 ^a$. It
is clear that $\xi _{kQ_1^c} (kQ_1^1) = (kQ_1^{1a})^*$, where $\xi
_{kQ_1^c}$ was defined in Lemma 1.7 of \cite{ZZC04}. Thus
$kQ_1^{1a}$ is a  Yetter-Drinfeld   $(kG)^*$-submodule of $kQ_1^a$, denoted by
$(kQ_1^{1a}, $ $ {\rm coad}(G, r, \overrightarrow {\rho}, u))$, which is
isomorphic to the dual of $(kQ_1^1, ad (G, r,$ $ \overrightarrow{
\rho}, u))$ as Yetter-Drinfeld $(kG)^*$-modules.

Therefore  we have the dual case of Theorem \ref {3}.

\begin{Proposition}\label{6.1}
Let $(G, r, \overrightarrow {\rho }, u)$ and $(G', r',
\overrightarrow{\rho'}, u')$ be two {\rm RSR}s. Then the following statements are equivalent:

{\rm(i)} ${\rm RSR} (G, r , \overrightarrow{\rho }, u)$ $\cong $
${\rm RSR} (G', r', \overrightarrow{\rho'}, u')$.

{\rm(ii)} There exists a Hopf algebra isomorphism $\phi:
(kG)^*\rightarrow (kG')^*$ such that  $(kQ_1^{1a},$ $ {\rm coad} (G,
r , \overrightarrow{ \rho}, u))\cong\ _{\phi}^{\phi^{-1}}(kQ_1'
{}^{1a}, {\rm coad} ( G', r', \overrightarrow{\rho'}, u'))$ as
Yetter-Drinfeld $(kG)^*$-modules.

{\rm(iii)} There is a Hopf algebra isomorphism $\phi:
(kG)^*\rightarrow (kG')^*$ such that ${\mathfrak B}(kQ_1^{1a}, $ $
{\rm coad}(G, r , \overrightarrow{ \rho}, u)$ $  \cong $ ${}_{\phi}
  ^{\phi ^{-1}} {\mathfrak B}  ( kQ_1' {}^{1a}, {\rm coad} ( G',
r', \overrightarrow{\rho'}, u'))$ as graded braided Hopf algebras in
$^{(kG)^*}_{(kG)^*} {\mathcal YD}$.
\end{Proposition}
\noindent {\bf Proof.} Obviously, {\rm (ii)} and Theorem \ref {3} {\rm (ii)} are equivalent;
{\rm (iii)} and Theorem \ref {3} {\rm (iii)} are equivalent. $\Box$

If $V = \oplus _{i=0}^\infty V_{(i)}$ is a graded vector space, and
${\rm dim} V_{(i)} < \infty$ for $0 \le i \le \infty$, then $V$ is
said to be locally finite. In this case, set $V^g =: \oplus
_{i=0}^\infty (V_{(i)})^* \subseteq V^*$, as in \cite {Sw69}.
\begin{Lemma}\label{6.2} Assume that  $H = \oplus _{n=0}^\infty  H_{(n)}$ is a locally finite
graded  Hopf algebra.

(i)  Then  $H^g =: \oplus _{n=0}^\infty (H_{(n)})^* \subseteq H^0$
is a locally finite graded  Hopf algebra.

(ii) If $G$ is a finite group and the coradical $H_0$ of $H$ is
$(kG)^*$, then the coradical $(H^g)_0 $ of $H^g$ satisfies $(H^g)_0
\cong kG.$
\end{Lemma}

\noindent {\bf Proof.} (i)  See Lemma 3.1.11 in \cite{Zh01}.

(ii)   Let $C_{n} =: \sum _{i =0}^n(H_{(i)})^*$, then $H^g= \sum
_{n=0} ^\infty C_n$ is a filtered coalgebra.  By Proposition
11.1.1 in \cite{Sw69}, $C_0 =((kG)^*)^* \cong kG$ contains the coradical
$(H^g)_0$ of  $H^g$. Consequently, $(H^g)_0 \cong  kG$. $\Box$

\begin{Proposition}\label{6.3} Let $G$ be a finite group.
Assume that   $V$ is a finite dimensional $(kG)^*$- ${\mathcal YD}$
module. Then

{\rm (i)} The coradical of $( \mathfrak B(V)\# (kG)^*)^g$ is
isomorphic to $kG$.

{\rm (ii)}  $( \mathfrak B(V)\# (kG)^*)^g$ is finite dimensional if
and only if $\mathfrak B(V)$ is finite dimensional.
\end{Proposition}
\noindent {\bf Proof.}  Let $R =:\mathfrak B(V) $ and the bosonization (or
biproduct) $R \#(kG)^* = \oplus _i^\infty (R_{(i)} \#(kG)^*)$ of $R$
and $(kG)^*$ is a local finite graded Hopf algebra with coradical
$(kG)^*$. Then $(R\# (kG)^*)^g =: \oplus _i^\infty (R_{(i)}
\#(kG)^*) ^*$ is a graded Hopf algebra  with   coradical $((R\#
(kG)^*)^g)_0 \cong kG$ by Lemma \ref {6.2}. This proves {\rm (i)}

For {\rm (ii)}, if  $(R\# (kG)^*)^g$ is finite dimensional, then  there exists a
natural number $m$ such that $(R\# (kG)^*)^g = \oplus _i^m (R_{(i)}
\#(kG)^*) ^*$. Consequently, $R \#(kG)^* = \oplus _i^m (R_{(i)}
\#(kG)^*)$ and $\mathfrak B(V)$ is finite dimensional. The converse is obvious. $\Box$
\vskip.1in
\noindent {\bf Remark.} This gives a method to decide whether Hopf algebra
with coradical  $(kG)^*$ is finite dimensional or not by means of pointed Hopf algebras. 


\section*{Acknowledgements}
We would like to thank Prof. N. Andruskiewitsch for his help.  The first and
third authors were financially supported by the Australian Research
Council, and the second author was supported by NSF of China
(10771183) and Sino-German project (GZ310). S.C.Z thanks the 
Department of Mathematics, University of Queensland for its hospitality.


\begin {thebibliography} {200}

\bibitem{Majid90} S Majid, {Quasi-triangular Hopf algebras and Yang-Bsxter equations},
Int. J. Mod. Phys. {\bf A 5} (1990), 1-91.



\bibitem{Ga} M.R. Gaberdiel, {An algebraic approach to logarithmic conformal field theory},
 Int. J. Mod. Phys. { \bf A 18} (2003), 4593-4638. 

\bibitem{AS98b} N. Andruskiewitsch and H. J. Schneider,
Lifting of quantum linear spaces and pointed Hopf algebras of order
$p^3$,  J. Algebra {\bf 209} (1998), 645--691.

\bibitem {AS02} N. Andruskiewitsch and H. J. Schneider, Pointed Hopf algebras, in {\em
New Directions in Hopf Algebras}, Math Sci. Res. Inst. Publ. {\bf 43} 
(Cambridge University Press, 2002), pp. 1-68. 


\bibitem{AS05} N. Andruskiewitsch and H. J. Schneider,
On the classification of finite-dimensional pointed Hopf algebras,
Ann. Math.  \textbf{171} (2010), 375--417 .


\bibitem{He04} I. Heckenberger, {\em The Weyl groupoid of a Nichols algebra
of diagonal type}, Invent. Math. \textbf{164} (2006), 175--188.


\bibitem{AZ07} N. Andruskiewitsch and S. Zhang, On pointed Hopf
algebras associated to some conjugacy classes in $S_n$, Proc. Amer. Math. Soc. {\bf 135} (2007), 2723-2731.

\bibitem{ZZC04} S. Zhang, Y.-Z. Zhang and H.X. Chen, Classification of PM quiver
Hopf algebras, J. Algebra Appl. {\bf 6} (2007), 919--950.

\bibitem{ZWW08} S. Zhang, M. Wu and H. Wang, Classification of ramification
systems for symmetric groups,  Acta Math. Sin. {\bf 51} (2008), 253--264.

\bibitem{DPR} R. Dijkgraaf, V. Pasquier and P. Roche,
Quasi Hopf algebras, group cohomology and orbifold models, Nucl. Phys. B Proc. Suppl. {\bf 18} (1991), 60--72.



\bibitem{CR02} C. Cibils and M. Rosso,  Hopf quivers, J. Algebra {\bf  254} (2002), 241-251.

\bibitem{CR97} C. Cibils and M. Rosso, Algebres des chemins quantiques,
Adv. Math. {\bf 125} (1997), 171--199.

\bibitem{ZZ07} S. Zhang and Y.-Z. Zhang, Structures and representations of
generalized path algebras, Algebr. Represent. Theor., {\bf 10} (2007), 117--134.

\bibitem{OZ04} F. Van Oystaeyen and P. Zhang, Quiver Hopf algebras, J. Algebra
{\bf 280} (2004), 577--589.

\bibitem{Ni78} W. Nichols, Bialgebras of type one, Comm. Algebra {\bf 6} (1978), 1521--1552.

\bibitem{Wo89} S.L. Woronowicz, Differential calculus on
 compact matrix pseudogroups (quantum groups), Comm. Math. Phys. {\bf 122} (1989), 125--170.

\bibitem{BD} Y. Bespalov and B. Drabant, Hopf (bi-)modules and crossed modules in
braided monoidal categories, J. Pure Appl. Algebra {\bf 123} (1998), 105--129.




\bibitem{Sw69} M. E. Sweedler, Hopf algebras, Benjamin, New York, 1969.

\bibitem{Zh01} S. Zhang, Braided Hopf algebras, Hunan Normal University Press, 1999; 
 The double bicrossproducts of braided Hopf algebras,  Comm.  Algebra {\bf  29} (2001),  31--66.
\end {thebibliography}
\end{document}